\documentclass[11pt]{amsart}
\usepackage{amssymb}

\theoremstyle{plain}
\newtheorem{theorem}{Theorem}

\newtheorem{example}{Example}

\theoremstyle{definition}

\theoremstyle{remark}

\numberwithin{equation}{section}

\begin{document} 

\title[Multi-variable Polynomial Solutions to Pell's Equation ]
       {Multi-variable Polynomial Solutions to Pell's Equation and 
Fundamental Units in Real Quadratic Fields}
\author{J. Mc Laughlin}
\address{Mathematics Department\\
       University of Illinois \\
        Champaign - Urbana, Illinois 61801}
\email{jgmclaug@math.uiuc.edu}
\date{December,11,2000}

\begin{abstract} 
Solving Pell's equation is of relevance in finding fundamental units in
real quadratic fields and for this reason polynomial solutions are 
of interest in that they can supply the fundamental units in infinite 
families of such fields.

In this paper an algorithm is described which allows one to construct, 
for each positive integer $n$, a finite collection, $\{F_{i}\}$,
of multi-variable 
polynomials (with integral coefficients), each satisfying a multi-variable 
polynomial Pell's equation
\[C_{i}^{2}-F_{i}H_{i}^{2}=(-1)^{n-1},
\]
where $C_{i}$ and $H_{i}$ are multi-variable polynomials with integral 
coefficients.
Each positive integer whose square-root has a regular continued fraction 
expansion with period $n+1$ lies in the range of one of these polynomials.
Moreover, the continued fraction expansion of these polynomials is given
explicitly as is the fundamental solution to the above
multi-variable polynomial Pell's equation.

Some implications for determining the fundamental unit in a wide class
of real quadratic fields is considered.

\end{abstract}

\maketitle
\section
{Introduction}\label{S:intro}
Solving Pell's equation is of relevance in finding fundamental units in
real quadratic fields and for this reason polynomial solutions are 
interesting in that they can supply the fundamental units in infinite 
families of such fields. 

There have been several papers written over the past thirty years which 
describe certain polynomials whose square roots have periodic 
continued fraction expansions which can be written down explicitly in terms 
of the coefficients and variables of the polynomials. See for example the 
papers of Bernstein~\cite{B76}, Levesque and Rhin~\cite{LR86},
 Madden~\cite{M01}, Van der Poorten~\cite{VDP94} and 
 Van der Poorten and Williams~\cite{VDPW99}.

In this paper an algorithm is described
which allows one to construct, for each positive integer $n$, a finite 
collection of multi-variable Fermat-Pell polynomials which have \emph{all}
positive integers whose square-roots have a continued fraction
expansion of period $n+1$ in their range. If 
$F_{i} := F_{i}(t_{0},t_{1},\cdots,t_{\lfloor \frac{n+1}{2} \rfloor})$ is any one of
these polynomials, the fundamental polynomial solution to the equation
\begin{equation}~\label{E:eq1}
C_{i}^{2} -F_{i}H_{i}^{2} = (-1)^{n-1}
\end{equation}
(where $C_{i}$ and $H_{i}$ are polynomials in the variables 
$t_{0},t_{1},\cdots,t_{\lfloor \frac{n+1}{2} \rfloor}$) can be found.
Moreover, the continued fraction expansion of $\sqrt{F_{i}}$ can be
written down when $t_{1},\cdots,
t_{\lfloor \frac{n+1}{2} \rfloor} \geq 0$  and 
$t_{0} > g_{i}(t_{1},\cdots,t_{\lfloor \frac{n+1}{2} \rfloor})$, a
certain rational function of these variables.
 Some implications
for single-variable Fermat-Pell polynomials are discussed as are  the 
implications for
writing down the fundamental units in a wide class of real quadratic
number fields.

\vspace{10pt}

\textbf{Definition:} 
 a multi-variable polynomial 
\[F:=F(t_{0},t_{1},\cdots, t_{k})\,\,
\in  \mathbb{Z}[\,t_{0},t_{1},\cdots, t_{k}],\,\,
k \geq 1 
\]
 is called
a \emph{multi-variable Fermat-Pell polynomial}
\footnote{These polynomials are called ``Fermat-Pell 
polynomials'' here to avoid confusion with 
``Pell Polynomials'' and also because Fermat investigated the ``Pell''
equation.} if there exists 
polynomials 
\[
C:=C(t_{0},t_{1},\cdots, t_{k}) \text{ and }
 H:=H(t_{0},t_{1},\cdots, t_{k}) \,\,
\in \mathbb{Z}[\,t_{0},t_{1},\cdots, t_{k}\,]
\]
such that either 
\begin{align}\label{E:pellm}
&C^{2} -F\,H^{2} =  1,\,\,\mbox{for all  $t_{i},\,\,\, 0 \leq
i \leq k$,\,\,\, or}\\
&C^{2} -F\,H^{2} = - 1,\,\,\mbox{for all $t_{i},\,\,\, 0 \leq
i \leq k$.} \notag
\end{align}
Such a triple of  polynomials 
$\{C,H,F\}$
satisfying equation~\eqref{E:pellm} constitute a  \emph{multi-variable
 polynomial solution} to Pell's equation.

\textbf{Definition:}
The multi-variable Fermat-Pell polynomial $F$ (as above)
 is said to have a 
\emph{multi-variable polynomial continued fraction expansion} if there
 exists a positive integer $n$, a real constant $T$, a rational 
function $g(t_{1},\cdots, t_{k})
\in \mathbb{Q}(t_{1},\cdots, t_{k})$  and 
polynomials 
$a_{0}:=a_{0}(t_{0},t_{1},\cdots, t_{k})\in
\mathbb{Z}[\,t_{0},t_{1},\cdots, t_{k}\,]$ and  
$a_{j}:=a_{j}(t_{1},\cdots, t_{k})\in
\mathbb{Z}[\,t_{1},\cdots, t_{k}\,],\,\,\,1\, \leq \, j\,\leq \, n$,  
which take only positive integral values for integral 
$t_{i} \geq T, 1\, \leq \, i\,\leq \, k$ and (possibly half-) integral  
$t_{0} > g(t_{1},\cdots, t_{k})$
such that 
\begin{equation*}
\sqrt{F} = [\,a_{0};
\overline{a_{1}, \cdots ,
a_{n},2a_{0}}], \text{ for all } t_{i}\,\text{'s in the ranges stated}
,\,0\leq i \leq k. 
\end{equation*}  
Remarks:\\
(1) From the point of view of simplicity it would be desirable to
replace the condition $t_{0} \geq g(t_{1},\cdots, t_{k})$ by $t_{0}
\geq T$ but it will be seen that for the polynomials examined here
that the former condition is more natural and indeed cannot be
replaced by the latter condition. \\
(2)The restriction that the $a_{i}(t_{1},\cdots, t_{k}) \geq 0$,
$1\, \leq \, i\,\leq \, n$
 may also
   seem artificial to some since negative terms can easily be removed
   from a continued fraction expansion (see, for example~\cite{VDP94})
but this changes the period of the continued fraction so  is avoided here.\\
(3) It may also seem artificial to have $a_{0}$ depend on a variable
    $t_{0}$ while the other $a_{i}$'s do not but this will also be
    seen to occur naturally.\\
(4)Finally, allowing $t_{0}$ to take half-integral values in some circumstances
   may also seem strange but this also will be seen to be natural and
indeed necessary.
\vspace{10pt}

\textbf{Definition:}
If, for all sets of integers $\{t_{0}',t_{1}',\cdots, t_{k}' \}$ 
satisfying $t_{0}' \geq g(t_{1}',\cdots, t_{k}')$ and 
 $t_{i}' \geq T$, $1\, \leq \, i\,\leq \, k$, 
\[
X = C_{i}(t_{0}',t_{1}',\cdots, t_{k}'),\,\,\,
Y = H_{i}(t_{0}',t_{1}',\cdots, t_{k}')
\]
 constitutes
the fundamental solution (in integers) to 
\[
X^2-F_{i}(t_{0}',t_{1}',\cdots, t_{k}')Y^2 = (-1)^{n-1}
\]
then $(C_{i}(t_{0},t_{1},\cdots, t_{k}),H_{i}(t_{0},t_{1},\cdots, t_{k}))$ is
termed the \emph{fundamental polynomial solution} to equation~\eqref{E:eq1}.

Standard notations are used: 
\[
a_{0}+\frac{1}{a_{1}+}\, \frac{1}{a_{2}+}\, 
\frac{1}{a_{3}+}\, \dots  \frac{1}{a_{N}}:= a_{0}+
\cfrac{1}{a_{1} + \cfrac{1}{a_{2} + \cfrac{1}{a_{3} +
\cdots   \cfrac{1}{a_{N} }}}}. 
  \]
To save space this continued fraction is usually written $[a_{0};a_{1},\cdots, a_{n}]$. 
The infinite periodic continued fraction with initial non-periodic part $a_{0}$ and
periodic part $a_{1} , ..., a_{n},2a_{0}$ is denoted by
$[\,a_{0};\overline{a_{1} , ......., a_{n},2a_{0}}]$. The $i$-th 
approximant of the continued fraction 
$[a_{0};a_{1},\cdots, ]$ is denoted by $P_{i}/Q_{i}$.

Repeated use will be made of  some basic
facts about  continued fractions, such as:
\begin{align}\label{E:bas1}
&P_{n}Q_{n-1}-P_{n-1}Q_{n} = (-1)^{n-1}, \\
&P_{n+1}= a_{n+1}P_{n}+P_{n-1}, \notag\\
&Q_{n+1}= a_{n+1}Q_{n}+Q_{n-1}, \notag
\end{align}
each of these relations being valid for $n=1,2,3 \cdots$.

Before coming to the main problem, it is necessary to first solve a 
related problem on  symmetric
strings of positive integers.

\section{A problem concerning Symmetric Sequences}

Question: For which symmetric sequences of positive integers $a_{1},...,a_{n}$ do there 
exist  positive integers $a_{0}$ and $D$ such that 
\begin{equation}\label{E:D}
\sqrt{D} = [\,a_{0};\overline{a_{1} , ......., a_{n},2a_{0}}]?
\end{equation}
Let $P_{i}/Q_{i}$ denote the $i$th approximant of the continued fraction 
\begin{equation}\label{E:incf}
\phantom{asd} 0+  \frac{1}{a_{1}+}\, \frac{1}{a_{2}+}\, 
\frac{1}{a_{3}+}\, \dots  \frac{1}{a_{n}}.
\end{equation}
 
By the well known correspondence between convergents and matrices
\begin{align*}
\left(
\begin{matrix}
0 & 1 \\
1     & 0
\end{matrix}
\right)
\left(
\begin{matrix}
a_{1} & 1 \\
1     & 0
\end{matrix}
\right)
\cdots
\cdots
\left(
\begin{matrix}
a_{n} & 1 \\
1     & 0
\end{matrix}
\right)
&=
\left(
\begin{matrix}
P_{n} & P_{n-1} \\
Q_{n} & Q_{n-1}
\end{matrix}
\right)\\
\Longrightarrow \left(
\begin{matrix}
0 & 1 \\
1     & 0
\end{matrix}
\right)
\left(
\begin{matrix}
a_{1} & 1 \\
1     & 0
\end{matrix}
\right)
\cdots
\cdots
\left(
\begin{matrix}
a_{n} & 1 \\
1     & 0
\end{matrix}
\right)
\left(
\begin{matrix}
0 & 1 \\
1     & 0
\end{matrix}
\right)
&=
\left(
\begin{matrix}
P_{n-1} & P_{n} \\
Q_{n-1} & Q_{n}
\end{matrix}
\right).\\
\end{align*}
Since the left side in the second equation is a symmetric sequence of 
symmetric matrices it follows that 
\begin{equation}\label{E:sym}
P_{n} = Q_{n-1}.
\end{equation}
Suppose $\sqrt{D}= [\,a_{0};\overline{a_{1} , ......., a_{n},2a_{0}}]
= a_{0}+\beta$, where  $\beta =$$[\,0;\overline{a_{1} , ......., a_{n},2a_{0}}] $ 
so that 
\begin{align*}
&\phantom{\Longrightarrow }\,\,\,\,
\beta =[\,0;a_{1} , ......., a_{n},2a_{0}+\beta] ,\\
&\Longrightarrow \beta = \frac{(2a_{0}+\beta)P_{n}+P_{n-1}}
			      {(2a_{0}+\beta)Q_{n}+Q_{n-1}} 
			= \frac{\beta P_{n} +(2a_{0}P_{n}+P_{n-1})}
				{\beta Q_{n} +(2a_{0}Q_{n}+Q_{n-1})},\\
&\Longrightarrow \beta^{2}Q_{n} + (2a_{0}Q_{n}+Q_{n-1} - P_{n})\beta -(2a_{0}P_{n}+P_{n-1}) 
				= 0,\\
&\Longrightarrow \beta^{2}Q_{n} + (2a_{0}Q_{n})\beta -(2a_{0}P_{n}+P_{n-1}) 
		= 0, \mbox{(by~\eqref{E:sym})} \\
&\Longrightarrow \sqrt{D}= a_{0} + \beta = \sqrt{a_{0}^{2} + \frac{2a_{0}P_{n}+P_{n-1}}{Q_{n}}}
\end{align*}
The problem now becomes one of determining for which 
symmetric sequences of positive integers $a_{1},...,a_{n}$ 
does there exist  positive 
integers $a_{0}$ such that 
$(2a_{0}P_{n}+P_{n-1})/Q_{n}$  is an integer.


\begin{theorem}\label{L:lem1}
There exists a positive integer $a_{0}$ such that $(2a_{0}P_{n}+P_{n-1})/Q_{n}$
is an integer if and only if $P_{n-1}Q_{n-1}$ is even.
\end{theorem}

\begin{proof}
$\Longleftarrow$ Suppose first of all that $P_{n-1}Q_{n-1}$ is even.
By equation~\eqref{E:bas1}
\[
P_{n}Q_{n-1}+(-1)^{n} = P_{n-1}Q_{n}
\]
(i) Suppose $n$ is even. Then 
$P_{n}Q_{n-1}P_{n-1}+P_{n-1} = P_{n-1}^{2}Q_{n}$.
Choose $t$ to be any integer or half-integer such that $tQ_{n}$ is an integer and 
$a_{0}:= Q_{n-1}P_{n-1}/2 + tQ_{n}> 0$. Then 
\[\frac{2a_{0}P_{n}+P_{n-1}}{Q_{n}} 
= \frac{Q_{n-1}P_{n-1}P_{n} + 2tP_{n}Q_{n} +P_{n-1}}{Q_{n}}
=2tP_{n} + P_{n-1}^{2}
\]
(ii)Similarly, in the case $n$ is odd, \,\, 
$-P_{n}Q_{n-1}P_{n-1}+P_{n-1} = -P_{n-1}^{2}Q_{n}$.
Choose  $t$ to be any integer or half-integer 
such that $tQ_{n}$ is an integer and 
$a_{0}:= -Q_{n-1}P_{n-1}/2 + tQ_{n}> 0$. In this case 
\[
\frac{2a_{0}P_{n}+P_{n-1}}{Q_{n}}=2tP_{n} - P_{n-1}^{2}
\]
$\Longrightarrow$ Suppose next that $P_{n-1}$ and $Q_{n-1}$ are both odd and that
there exists a positive integer $a_{0}$ such that $(2a_{0}P_{n}+P_{n-1})/Q_{n}$  
is a positive integer, $m$, say. Using  ~\eqref{E:bas1} and ~\eqref{E:sym} it follows that $Q_{n}$ is even. Then 
 $2a_{0}P_{n}+P_{n-1} = mQ_{n}$ implies $P_{n-1}$ is even - a contradiction.
\end{proof}
Remarks:\\
(i)Note that this process gives all $a_{0}$ such that 
$(2a_{0}P_{n}+P_{n-1})/Q_{n}$  is an integer. Indeed,
\begin{align*}
&(2a_{0}P_{n}+P_{n-1})/Q_{n} = k,\mbox{ an integer}\\
&\Longleftrightarrow 2a_{0}P_{n}Q_{n-1} = -P_{n-1}Q_{n-1} + kQ_{n}Q_{n-1},\\
&\Longleftrightarrow 2a_{0}(-1)^{n-1}=
2a_{0}(P_{n}Q_{n-1}- P_{n-1}Q_{n}),\,\, \mbox{( by~\eqref{E:bas1})} \\
&\phantom{\Longleftrightarrow 2a_{0}(-1)^{n-1}}
 = -P_{n-1}Q_{n-1} + Q_{n}(kQ_{n-1}-2a_{0}P_{n-1}),\\
&\Longleftrightarrow a_{0} = (-1)^{n-1}\left(\frac{-P_{n-1}Q_{n-1}}{2}
			+Q_{n}\frac{kQ_{n-1}-2a_{0}P_{n-1}}{2}\right).
\end{align*}
Notice also that if there is one such $a_{0}$ that there are
infinitely many of them.

(ii)Notice that, with $P_{n},P_{n-1},Q_{n}$ and $Q_{n-1}$  as defined
above, if there exists a positive integer $D$ satisfying~\eqref{E:D} then
$D = p(t_{0})$, for some allowed $t_{0}$, where 
\begin{align*}
&p(t) = \left(\frac{Q_{n-1}P_{n-1}}{2} + tQ_{n}\right)^{2} + 2tP_{n} + P_{n-1}^{2},
	\,\,	t> \frac{-Q_{n-1}P_{n-1}}{2Q_{n}},\,\,\,\, \mbox{($n$ even),} \notag \\
&p(t) = \left(\frac{-Q_{n-1}P_{n-1}}{2} + tQ_{n}\right)^{2} + 2tP_{n} - P_{n-1}^{2},
	\,\,	t> \frac{Q_{n-1}P_{n-1}}{2Q_{n}},\,\,\, \mbox{($n$ odd).} \notag
\end{align*}

The above theorem suggests a simple algorithm for deciding if, for a
given symmetric
sequence of positive integers $a_{1}, \cdots, a_{n}$, there exist positive
integers $a_{0}$ and $D$ such that~\eqref{E:D} holds. Notice that all that 
matters is the parity of the $a_{i}$ so all calculations can be done in 
$\mathbb{Z}_{2}$. First of all define the following matrices: 
\[
J=\left(
\begin{matrix}
0 & 1\\
1 & 0
\end{matrix}
\right),\,\,
K=\left(
\begin{matrix}
1 & 1\\
1 & 0
\end{matrix}
\right)\,\,
\mbox{and} \,\, I = 
\left(
\begin{matrix}
1 & 0\\
0 & 1
\end{matrix}
\right).
\]
Convert the sequence $a_{1},a_{2},\cdots,a_{n}$ to a sequence of $J$-
and $K$-matrices,
according to whether each $a_{i}$ is odd (replace by a $K$) or even (replace by a $J$).
Prefix a $J$-matrix (to account for the initial $0$ in the 
continued fraction~\eqref{E:incf}). Multiply this sequence together (modulo $2$) 
using the facts that
$J^{2} = K^{3}=I, \mbox{ and } J K = K^2 J$. 

The final matrix $ \equiv  \left(
\begin{matrix}
* & 1\\
* & 1
\end{matrix}
\right)\mod{2} 
\Longleftrightarrow$ there do not exist positive integers $a_{0}$ and $D$ such 
		that~\eqref{E:D} holds. 

\begin{example}\label{E:EX}
 Do there exist positive integers $a_{0}$ and $D$  such that 
\[ \sqrt{D}= [\,a_{0};\overline{22,34,97,32,15,17,17,15,32,97,34,22,2a_{0}}]?
 \]
\end{example}
As described above convert the sequence
$22,34,97,32,15,17,17,15,32,97$,
$34,22$ to
a sequence of $J$- and $K$-matrices, prefix a $J$-matrix and multiply
the sequence together:
{\allowdisplaybreaks
\begin{align*}
\underbrace{JJ} JKJ \underbrace{KKK} KJK \underbrace{JJ}=
JK (JK) JK &= J \underbrace{K(K^{2}} \underbrace{J)J} K\\
&= JK =\left(
\begin{matrix}
1 & 0\\
1 & 1
\end{matrix}
\right).
\end{align*}
}
Therefore there do exist positive integers $a_{0}$ and $D$ 
such that 
\[ \sqrt{D}= [\,a_{0};\overline{22,34,97,32,15,17,17,15,32,97,34,22,2a_{0}}].
 \]

\section{Multi-variable Fermat-Pell Polynomials}
\textbf{Definition:} If $\{\,a_{1},\cdots,a_{n}\}$ is a symmetric zero-one
sequence such that 
\begin{equation*}
\left(
\begin{matrix}
0 & 1\\
1 & 0
\end{matrix}
\right)
\prod_{i=1}^{n}
\left(
\begin{matrix}
a_{i} & 1\\
1 & 0
\end{matrix}
\right)
\not \equiv
\left(
\begin{matrix}
* & 1\\
* & 1
\end{matrix}
\right) \mod{2}
\end{equation*} 
then the sequence $\{\,a_{1},\cdots,a_{n}\}$ is termed a
\emph{permissible} sequence. Let $r(n)$ denote the number of
permissible sequences of length $n$.

Note: It is not difficult to show that 
$r(2m)= ((-1)^{m} + 2^{m + 1})/3$ and that  $r(2m+1)= 
( (-1)^{m} + 5\times2^{m})/3$.

  If $D$ is a positive integer
such that $\sqrt{D} = [\,a_{0};\overline{a_{1} , .......,
a_{n},2a_{0}}]$ then $\{\,a_{1},\cdots,a_{n}\} \mod{2}$ must equal one
of the above permissible sequences and $D$ is said to be
\emph{associated} with this permissible sequence . The collection of
all positive integers
associated with a particular permissible sequence is termed the 
\emph{parity class} of this permissible sequence. Sometimes, if there
is no danger of ambiguity, these collections of positive integers will
be referred to simply as \emph{parity classes}.


\begin{theorem}\label{E:1par}
(i)For  each positive integer $n$ there exists a finite collection of 
multi-variable Fermat-Pell polynomials 
$F_{j}(t_{0},t_{1},\cdots,t_{\lfloor \frac{n+1}{2} \rfloor}), 1 \leq j
\leq r(n)$,
such that each positive integer whose square root has a continued
fraction expansion with period $n+1$ lies in the range of exactly one
of these polynomials.
Moreover, these polynomials can be constructed;

(ii) These polynomials have a polynomial continued
fraction expansion which can be explicitly determined;

(iii) The fundamental polynomial solution \\
$C=C_{j}(t_{0},t_{1},\cdots,t_{\lfloor \frac{n+1}{2} \rfloor}), \,\,\,
H=H_{j}(t_{0},t_{1},\cdots,t_{\lfloor \frac{n+1}{2} \rfloor})$ 
to
\begin{equation}\label{E:fund}
C^{2} -
F_{j}(t_{0},t_{1},\cdots,t_{\lfloor \frac{n+1}{2} \rfloor})
H^{2} = (-1)^{n-1}
\end{equation}
exists and can be explicitly determined.
\end{theorem}

\begin{proof}
(i)The proof will be by construction.

\emph{Step 1}: Find all permissible  sequences. 
This will involve checking $2^{\lfloor \frac{n+1}{2} \rfloor}$ zero-one sequences
in a way similar to the example~\eqref{E:EX} above.

\emph{Step 2}: 
For each permissible sequence 
$\{a_{1},\cdots, a_{n}\}$
create a new symmetric polynomial sequence 
$\{a_{1}(t_{1}),a_{2}(t_{2}),\cdots,a_{n-1}(t_{2}), a_{n}(t_{1})\}$
by replacing each $a_{i}$ and its partner $a_{n+1-i}$ in the symmetric 
sequence 
 by $a_{i}(t_{i}) = a_{n+1-i}(t_{i})= 2 t_{i}+1$ if $a_{i}=1$  and by 
$a_{i}(t_{i}) = a_{n+1-i}(t_{i})= 2 t_{i}+2$ if $a_{i}=0$.
This new sequence will sometimes be referred to as the sequence 
$\{a_{1},\cdots, a_{n}\}$,
if there is no danger of ambiguity.
Each of the integer variables  $t_{i}$ (in the polynomial being constructed)
will be allowed to
vary independently over the range $0 \leq t_{i} < \infty$
and each of the new $a_{i}$'s will keep the same parity and stay positive.

\emph{Step 3}
As in~\eqref{E:incf}, form the continued fraction
\[
\phantom{asd} 0+  \frac{1}{a_{1}(t_{1})+}\, \frac{1}{a_{2}(t_{2})+},\,
 \dots,  
\frac{1}{a_{n-1}(t_{2})+}\,\frac{1}{a_{n}(t_{1})}
\]
and calculate $P_{n},Q_{n},P_{n-1}$ and $Q_{n-1}$ for this  polynomial continued 
fraction, where these expressions are now polynomials in the $t_{i}$'s. 

\emph{Step 4} 
Construct $F_{j}:=
F_{j}(t_{0},t_{1},\cdots,t_{\lfloor \frac{n+1}{2} \rfloor})$, the multi-variable
Fermat-Pell polynomial corresponding to the particular parity sequence under 
consideration.
This is simply done by defining
\begin{align}\label{E:tot}
F_{j}:=
\begin{cases}
\left(\frac{Q_{n-1}P_{n-1}}{2} + t_{0}Q_{n}\right)^{2} + 2t_{0}P_{n} + P_{n-1}^{2},
\,\, \mbox{($n$ even)}\\
\left(\frac{-Q_{n-1}P_{n-1}}{2} +t_{0}Q_{n}\right)^{2}+2t_{0}P_{n}-P_{n-1}^{2},
\,\, \mbox{($n$ odd)}
\end{cases}
\end{align}
where $(-1)^{n+1}Q_{n-1}P_{n-1}/(2Q_{n}) \,\,< t_{0} \,\,<\, \infty$ 
and $t_{0}$ can take half-integral
 values if $Q_{n}$ is even and otherwise takes integral values.

Every positive integer whose square root has a continued fraction
 expansion with period $n+1$ lies in the range of exactly one of these
 polynomials. That these polynomials are multi-variable Fermat-Pell
 polynomials follows from equation~\eqref{E:pelfun} below.

(ii) With $t_{0}$ in the range given, then 
\begin{equation*}
\sqrt{F_{j}} = [\,a_{0}(t_{0},t_{1},\cdots,t_{\lfloor \frac{n+1}{2} \rfloor});
\overline{a_{1}(t_{1}) , \cdots, 
a_{n}(t_{1}),
2a_{0}(t_{0},t_{1},\cdots,t_{\lfloor \frac{n+1}{2} \rfloor})}],
\end{equation*}
for all $t_{i} \geq 0$. Here 
\begin{align}\label{E:a0}
a_{0}=a_{0}(t_{0},t_{1},\cdots,t_{\lfloor \frac{n+1}{2} \rfloor}):=
\begin{cases}
\frac{Q_{n-1}P_{n-1}}{2} + t_{0}Q_{n},
\,\, \mbox{($n$ even)}\\
\frac{-Q_{n-1}P_{n-1}}{2}+ t_{0}Q_{n},
\,\, \mbox{($n$ odd)}.
\end{cases}
\end{align}

(iii) Notice  (using~\eqref{E:bas1} and~\eqref{E:sym}) that 
\begin{equation}\label{E:pelfun}
(a_{0}Q_{n}+P_{n})^{2}-(a_{0}^{2} + (2a_{0}P_{n}+P_{n-1})/Q_{n})Q_{n}^{2}= (-1)^{n-1}.
\end{equation} 
To see that
$(a_{0}Q_{n}+P_{n},Q_{n})$ is the fundamental solution
to~\eqref{E:fund}, notice that 
\[
\sqrt{F_{j}} = [\,a_{0}(t_{0},t_{1},\cdots,t_{\lfloor \frac{n+1}{2} \rfloor});
\overline{a_{1}(t_{1}) , \cdots, 
a_{n}(t_{1}),
2a_{0}(t_{0},t_{1},\cdots,t_{\lfloor \frac{n+1}{2} \rfloor})}].
\]
This has 
period $n+1$ and the $n$th approximant is 
$a_{0} + P_{n}/Q_{n} = (a_{0}Q_{n}+P_{n})/Q_{n}$
and by the theory of the Pell equation 
$(a_{0}Q_{n}+P_{n},Q_{n})$ is the fundamental solution
to~\eqref{E:fund}.
\end{proof}

As regards fundamental units in quadratic fields there is the
following theorem on page 119 of~\cite{N90}:
\begin{theorem}
Let $D$ be a square-free, positive rational integer and let 
$K=\mathbb{Q}(\sqrt{D})$. Denote by $\epsilon_{0}$ the fundamental
unit of $K$ which exceeds unity, by $s$ the period of the continued
fraction expansion for $\sqrt{D}$, and by $P/Q$ the ($s-1$)-th
approximant of it.

If $D \not \equiv 1 \mod{4}$ or $D \equiv 1 \mod{8}$, then 
\[\epsilon_{0} = P + Q \sqrt{D}.
\] 
However, if $D \equiv 5 \mod{8}$, then 
\[\epsilon_{0} = P + Q \sqrt{D}.
\] 
or 
\[\epsilon_{0}^{3} = P + Q \sqrt{D}.
\] 
Finally, the norm of $ \epsilon_{0}$ is positive if the period $s$ is 
even and negative otherwise.
\end{theorem}
It is easy, working modulo $4$, to determine simple conditions 
(on $t_{0}$ ) which make 
$F_{j} \equiv 2
\,\,
\mbox{ or}\,\, 3\,\, \mod{4}$   
and thus to say further, for a particular set of choices of 
$t_{1},\cdots,t_{\lfloor \frac{n+1}{2} \rfloor}$ and for all odd or
even $t_{0}$, that if 
$F_{j}$ is
square-free, then 
$a_{0}Q_{n}+P_{n}+ \sqrt{
F_{j}}Q_{n}$  
 is the fundamental unit 
in $ \mathbb{Q}[\,\sqrt{F_{j}}\,\,]$.
 For example, suppose that $n$ is even and that the original $Q_{n-1}$
determined from the permissible zero-one sequence is also even 
(so that $P_{n-1}$ and $Q_{n}$ are both
odd and $P_{n}=Q_{n-1}$ is even). 
Then the multi-variable form of $Q_{n-1}$
evaluated in \emph{Step 3} will also have all even coefficients. Suppose
$\frac{Q_{n-1}}{2} \equiv c_{0} + \sum t_{i^{'}}\mod{2} $. 
(Here $c_{0}$ may be $0$
and the sum $\sum t_{i^{'}}$ may contain some, all or none of the
$t_{i}$'s )  It is easy to see that
$F_{j}(t_{0},t_{1},\cdots,t_{\lfloor \frac{n+1}{2} \rfloor}) \equiv 
(c_{0} + \sum t_{i^{'}} + t_{0})^{2}+1 \mod{4}$. Even more simply, if the
original original $Q_{n-1}$ as in \emph{Step 1} is odd (here also the case $n$ is
even is considered) then $P_{n-1}$ as evaluated in \emph{Step 3} is even and 
 it is not difficult to show that in fact $ P_{n-1} \equiv 2
\mod{4}$ (since for $n$ even $P_{n}Q_{n-1}-P_{n-1}Q_{n} = -1$) and
that $Q_{n}$ is odd,
which leads to 
$F_{j}(t_{0},t_{1},\cdots,t_{\lfloor \frac{n+1}{2} \rfloor}) \equiv
t_{0}^{2}+1 \mod{4}$.
Similar relations hold in the case where $n$ is odd.

 The polynomials constructed in theorem~\eqref{E:1par} take values
in only one parity class, if all the variables are positive. However,
given any two parity classes,
there are multi-variable Fermat-Pell polynomials that take values in
those two classes.


\begin{theorem}\label{E:2par} 
Let $n$ be any fixed positive integer large enough so that the set of 
positive integers whose square roots have a continued fraction 
expansion of period $n+1$ can be divided into more than one parity class.

(i) Given any two parity classes of integers whose square roots
have continued fraction expansions of period $n+1$, 
there are multi-variable Fermat~\--Pell polynomials, which can be
constructed, that take values in
both parity classes;

(ii) These polynomials have a polynomial continued
fraction expansion which can be explicitly determined;

(iii) If $F=F(t_{0},c,t_{1},\cdots, t_{\lfloor \frac{n+1}{2} \rfloor})$ is
any such polynomial then the 
 fundamental polynomial solution 
\[
C=C(t_{0},c,t_{1},\cdots, t_{\lfloor \frac{n+1}{2} \rfloor}), 
H=H(t_{0},c,t_{1},\cdots, t_{\lfloor \frac{n+1}{2} \rfloor})
\]
to
\begin{equation}
C^{2} -FH^{2} = (-1)^{n-1}
\end{equation}
can be explicitly determined.
\end{theorem}

\begin{proof}
As in \emph{Step 2} in theorem~\eqref{E:1par}  a polynomial sequence 
$\{a_{1},\cdots, a_{n}\}$ is created. 
Suppose $L_{1}=\{b_{1} \cdots,b_{n}\}$ and $L_{2}=\{c_{1},\cdots,
c_{n}\}$ 
are the permissible sequences associated
with the two parity classes.
Let $i_{1},\cdots,i_{k}$ be those positions 
$\leq \lfloor \frac{n+1}{2} \rfloor$  at which the sequences
 agree. For each of these $i_{r}$'s set 
$a_{i_{r}}(t_{i_{r}}) = a_{n+1-i_{r}}(t_{i_{r}}) = 
2 t_{i_{r}}+1$, if $c_{i_{r}}$ is odd and set  
$a_{i_{r}}(t_{i_{r}}) =a_{n+1-i_{r}}(t_{i_{r}}) =  2 t_{i_{r}}+2$, 
if $c_{i_{r}}$ is even. Subdivide the remaining positions
(those positions $\leq \lfloor \frac{n+1}{2} \rfloor$ at which $L_{1}$
and $L_{2}$ 
differ) into two 
subsets: those at which $L_{1}$ has a $0$ and
$L_{2}$ has a $1$ and those at which $L_{1}$ has a $1$ and
$L_{2}$ has a $0$. 

Suppose $i_{j}$ is a position of the first kind. Let 
$a_{i_{j}}(c,t_{i_{j}}) = a_{n+1-i_{j}}(c,t_{i_{j}})$
$ = c+2+2t_{i_{j}}$. 
Repeat this for all the positions $i_{j}$ in this first set. 
 Likewise,  Suppose $i_{j}$ is a position of the second kind.
In this case let  
$a_{i_{j}}(c,t_{i_{j}}) = a_{n+1-i_{j}}(c,t_{\lfloor \frac{n+1}{2} \rfloor})
$ $= c+1+2t_{i_{j}}$.This is also repeated for all the positions $i_{j}$ in this second set.
\emph{Step 3} and \emph{Step 4} are then carried out as above. 
The rest of the proof is  identical to theorem~\eqref{E:1par}. 
Denote the polynomial produced by 
\begin{equation}\label{E:totb}
F:=F(t_{0},c,t_{1},\cdots, t_{\lfloor \frac{n+1}{2} \rfloor}).
\end{equation}
As in theorem~\eqref{E:1par}, if $c$ and all the $t_{i}$'s are
non-negative, $1 \leq i \leq \lfloor \frac{n+1}{2} \rfloor$
 and  $t_{0}>(-1)^{n+1}Q_{n-1}P_{n-1}/(2Q_{n})$ then 
\[
\sqrt{F}= [\,a_{0};\overline{a_{1} , ......., a_{n},2a_{0}}],
\]
where the $a_{i}$'s, \,\,$1 \leq i \leq n$ are as defined just above
and $a_{0}$ is as defined in equation~\eqref{E:a0}.

 Under these conditions also
 the parity class of $F(t_{0},c,t_{1},\cdots, t_{\lfloor
\frac{n+1}{2} \rfloor})$
will depend only on the parity of $c$.
As in theorem ~\eqref{E:1par} the fundamental polynomial solution to 
\[
C^{2} -F(t_{0},c,t_{1},\cdots, t_{\lfloor \frac{n+1}{2} \rfloor})H^{2} =
(-1)^{n-1}
\]
is given by $C=a_{0}Q_{n}+P_{n},\,\,\,H=Q_{n}$.
\end{proof}

\section{A Worked Example}

As an example, consider those positive integers whose square-roots have  
continued fraction expansion
with period of length 9. Thus the symmetric  part of the period has length $8$
and it is necessary to check the $2^{4}=16$ zero-one sequences to determine
which are permissible. (This checking is done in essentially the same way as 
in Example~\ref{E:EX} above.)
There are $11$ valid sequences:
{\allowdisplaybreaks
\begin{align*}
&\phantom{asdf}\\
&{0, 0, 0, 0, 0, 0, 0, 0}\notag	\\
&{0, 0, 0, 1, 1, 0, 0, 0}\notag\\
&{0, 0, 1, 1, 1, 1, 0, 0}\notag\\
&{0, 1, 0, 0, 0, 0, 1, 0}\notag\\
&{0, 1, 0, 1, 1, 0, 1, 0}\notag\\
&{0, 1, 1, 1, 1, 1, 1, 0}\notag\\
&{1, 0, 0, 1, 1, 0, 0, 1}\notag\\
&{1, 0, 1, 0, 0, 1, 0, 1}\notag\\
&{1, 0, 1, 1, 1, 1, 0, 1}\notag\\
&{1, 1, 0, 0, 0, 0, 1, 1}\notag\\
&{1, 1, 1, 0, 0, 1, 1, 1}\notag
\end{align*}
}
The ninth of these is considered in more detail
(Each of the others can be dealt with in a similar way). 
For clarity the letters $a,b,c$ and $d$ are used instead of $t_{1},t_{2},t_{3}$
and $t_{4}$. 
Evaluating the continued fraction
\begin{equation}\label{E:ex}
 0+  
\frac{1}{2a+1+}\, \frac{1}{2b+2+}\, 
\frac{1}{2c+1+}\,\frac{1}{2d+1+}\,
\frac{1}{2d+1+}\, \frac{1}{2c+1+}\, 
\frac{1}{2b+2+}\,   \frac{1}{2a+1}
\end{equation}
 it is found that 
{\allowdisplaybreaks
 \begin{align*}
 P_{8}&=Q_{7}= -1 - 2\,d +\\
& 2\,( 3 + 4\,a + 2\,b + 4\,a\,b ) 
   ( 4 + 3\,b + 4\,c + 4\,b\,c + 6\,d + 4\,b\,d + 
     8\,c\,d + 8\,b\,c\,d )  \\&+
  4\,( 3 + 2\,b + 4\,( 1 + b ) \, c )
     \,( 2 + b + 3\,c + 2\,b\,c + 
     a\,( 3 + 2\,b + 4\,( 1 + b ) \,c
        )  ) \times \\
&\phantom{dfdsfsdafsafsafdsgggdfgdfgdfgdssfdsgfsddfgsdfgd}
( 1 + 2\,d + 2\,d^2 ),\\
 P_{7}&=4\,( 1 + b ) \,
   ( 4 + 3\,b + 4\,c + 4\,b\,c + 6\,d + 4\,b\,d + 
     8\,c\,d + 8\,b\,c\,d )  + \\
 &\phantom{dfdsfsdafsafsafdsgggdfg}
 2\,{( 3 + 2\,b + 4\,( 1 + b ) \,c ) 
       }^2\,( 1 + 2\,d + 2\,d^2 )\text{ and}\\
 Q_{8}&=8\,( 2 + b + 3\,c + 2\,b\,c + 
       a\,( 3 + 2\,b + 4\,( 1 + b ) \,c
          )  )^2\,
   ( 1 + 2\,d + 2\,d^2 )  + \\
  &( 3 + 4\,a + 2\,b + 4\,a\,b ) \,
   ( 3 + 4\,c + 4\,d + 8\,c\,d + 
     ( 2 + 4\,a ) \,
      ( 4 + 3\,b +\\
&\phantom{dfdsfsdafsafsafdsgggdfgdfgj} 4\,c + 4\,b\,c + 6\,d + 4\,b\,d + 
        8\,c\,d + 8\,b\,c\,d )  ).
\end{align*}
}
Since $n$ is $8$ (even) and $Q_{8}$ is odd (so $t_{0}$ 
cannot take half-integer values),  in this case 
$F_{9}(t_{0},a,b,c,d)$ is defined by 
\begin{equation}\label{E:expoly}
F_{9}(t_{0},a,b,c,d)=(Q_{7}P_{7}/2 + t_{0}Q_{8})^2 + 2t_{0}P_{8} + P_{7}^{2}
\end{equation} 
and 
\begin{multline*}
\sqrt{F_{9}(t_{0},a,b,c,d)}= 
[Q_{7}P_{7}/2 + t_{0} Q_{8};\overline{
2a+1 , 2b+2, 
2c+1 ,  2d+1 ,
2d+1,}\\
 \overline{2c+1,
2b+2 , 2a+1 , 2(Q_{7}P_{7}/2 + t_{0} Q_{8})}],
\end{multline*}
 this expansion being valid for all $a,b,c,d \geq 0$ and all 
$t_{0}>-Q_{7}P_{7}/(2Q_{8})$ and in particular for all $t_{0} \geq 0$.
In these ranges 
\[
C=  (Q_{7}P_{7}/2 + t_{0}Q_{8})Q_{8}+P_{8},\,\,H=Q_{8}
\]
 gives the fundamental polynomial 
solution to
\[
C^{2}-F_{9}H^{2} = -1.
\]
\[
F_{9}(t_{0},a,b,c,d)=(Q_{7}P_{7}/2 + t_{0}Q_{8})^2 + 2t_{0}P_{8} + P_{7}^{2} 
\equiv (1+t_{0}^{2})\mod{4}
\]
 so that if   
$(Q_{7}P_{7}/2 + t_{0}Q_{8})^2 + 2t_{0}P_{8} + P_{7}^{2}$ is a square-free
number for some particular $a,b,c,d\, \geq \,0$ and some odd 
$t_{0}>-Q_{7}P_{7}/(2Q_{8}) $, then  
\[
(Q_{7}P_{7}/2 + t_{0}Q_{8})Q_{8}+P_{8} + 
\sqrt{(Q_{7}P_{7}/2 + t_{0}Q_{8})^2 + 2t_{0}P_{8} + P_{7}^{2}}Q_{8}
\]
 is the fundamental unit 
in $\mathbb{Q}\left(\,\sqrt{(Q_{7}P_{7}/2 + t_{0}Q_{8})^2 + 2t_{0}P_{8} + P_{7}^{2} }\right)$.

\section{Mystification, Fermat-Pell polynomials of a single variable and more on 
odd-even}

Clearly it is possible to ``mystify'' this process by replacing each $t_{i}$
by some polynomial $g_{i}(t_{i})$ taking only positive values or by replacing $2t_{i}$
(recalling that the continued fraction expansion contains only terms like
$2t_{i} +1$ or $2t_{i} +2$)
 by some polynomial  $g_{i}(t_{i})$ taking only 
even non-negative values or by setting $t_{i}=
t_{i}(X_{1},X_{2}, \cdots, X_{k}),\,\, 
1 \leq i \leq \lfloor \frac{n+1}{2} \rfloor $,  a 
polynomial in the $X_{j}$'s taking only positive values,
where  the $X_{j}$'s can be independent variables and $k$ can be as large as desired
and so on.

Finally of course one can obtain single-variable Fermat-Pell 
polynomials by replacing the original variables 
$t_{0}, t_{i},1 \leq i \leq \lfloor \frac{n+1}{2} \rfloor $ by 
polynomials in a single variable. If  it is desired that 
the period of the continued fraction expansion of the new single-variable 
Fermat-Pell polynomial  should stay the same as that of 
the originating multi-variable polynomial  then the domain of the single
variable should be  restricted so that  the polynomials
replacing each of the $t_{i}$'s take only 
positive  values as in the multi-variable case and the polynomial 
replacing $t_{0}$
must be such that the $a_{0}$ term stays positive for all allowed
values of the new single variable.
 
For example, letting $a=s,b=0,c=s,d=0$ and $t_{0}=s$ in the
 polynomial~\eqref{E:expoly}
 above
produces the single-variable Fermat-Pell polynomial 
{\allowdisplaybreaks
\begin{multline*}
g(s)=639557 + 6858268\,s + 33078145\,s^2 +\\ 94534688\,s^3 + 
  177380352\,s^4 +
 228442240\,s^5 + 204593408\,s^6 + \\
  125870080\,s^7 + 50925568\,s^8 + 12238848\,s^9 + 
  1327104\,s^{10}
\end{multline*}
}
which has the continued fraction expansion (valid for all $s \geq 0$)
\begin{multline*}
\sqrt{g(s)}= 
[799 + 4289\,s + 9184\,s^2 + 9856\,s^3 + 5312\,s^4 + 
  1152\,s^5;\\
\overline{
2s+1 , 2, 
2s+1 ,  1 ,
1 , 2s+1 , 
2 , 2s+1 ,}\phantom{asdffsdafs}\\
\overline{
2(799 + 4289\,s + 9184\,s^2 + 9856\,s^3 + 5312\,s^4 + 
  1152\,s^5 )}].
\end{multline*}
$g(s) \equiv (1+ s^2) \mod{4}$ so when $s$ is odd and positive
 and $g(s)$ is square-free 
\begin{multline*}
51982 + 534625\,s + 2429840\,s^2 + 6408000\,s^3 + \\
  10812928\,s^4 + 12115200\,s^5 +
 9019392\,s^6 + 
  4304896\,s^7 + 1196032\,s^8 +\\ 147456\,s^9 
+
\sqrt{g(s)}(65 + 320\,s + 576\,s^2 + 448\,s^3 + 128\,s^4)
\end{multline*}
is the fundamental unit in $\mathbb{Q}[\,\sqrt{g(s)}]$.
For example, letting $s=1$ gives that 
$47020351+ 1537\sqrt{935888258}$ is the fundamental unit in 
$\mathbb{Q}[\,\sqrt{935888258}]$.

Starting with the continued fraction
\begin{equation*}
 0+  
\frac{1}{2a+1+}\, \frac{1}{2b+2+}\, 
\frac{1}{c+2e+}\,   \frac{1}{2d+1+}\,
\frac{1}{2d+1+}\, \frac{1}{c+2e+}\, 
\frac{1}{2b+2+}\,   \frac{1}{2a+1}\,
\end{equation*}
and following the same steps as above with the continued fraction~\eqref{E:ex}
 a multi-variable Fermat-Pell polynomial is developed which takes
 values in the parity classes associated with permissible sequences 7 and 9. Letting
$a=b=d=e=t=0$ one gets the single-variable Fermat-Pell polynomial
\begin{multline*}
g(c)=4325 + 28140\,c + 83652\,c^{2} + 147440\,c^{3} + 168000\,c^{4} + \\
  126528\,c^{5} + 61504\,c^{6} + 17664\,c^{7} + 2304\,c^{8}
\end{multline*}
with continued fraction expansion
\begin{multline*}
\sqrt{g(c)}=[ 65 + 214 c + 288 c^{2} + 184 c^{3} + 48 c^{4}; \\
\overline{1,2,c,1,1,c,2,1,2(65 + 214 c + 288 c^{2} + 184 c^{3} + 48 c^{4})}],
\end{multline*}
valid for $c \,\geq \, 1$.

\section{Concluding Remarks}
Every Fermat-Pell polynomial in one variable, $s$ say, 
that eventually has a continued fraction expansion
of fixed period length  can be found from~\eqref{E:tot}, if it takes
values in only one parity class for all sufficiently large $s$, and 
from \eqref{E:totb}, if it takes
values in two parity class for all sufficiently large $s$.
(Recall remark (i) after theorem~\eqref{L:lem1})

Of course none of this does anything to answer Schinzel's question of 
whether every Fermat-Pell polynomial in one variable has a continued fraction 
expansion. Neither does it provide a criterion (such as Schinzel's in
the degree-two case) for deciding if a polynomial of arbitrarily high
even degree is a Fermat-Pell polynomial.
 Perhaps it raises another question - Does every multi-variable 
Fermat-Pell polynomial have a continued fraction expansion? Does every multi-variable 
Fermat-Pell polynomial have a continued fraction expansion, assuming every 
Fermat-Pell polynomial in one variable does?


\begin{thebibliography}{9}
\bibitem{B76}
	Bernstein, Leon.
	\emph{ Fundamental units and cycles in the period of real 
	quadratic number fields. I.} Pacific J. Math. \textbf{63} (1976), no. 1, 37--61. 

\bibitem{B76a}
Bernstein, Leon.
\emph{ Fundamental units and cycles in the period of real 
quadratic number fields. II.} Pacific J. Math. \textbf{63} (1976), no. 1, 63--78. 

\bibitem{E48}
	L. Euler, (Translated by John D. Blanton)
	  \emph{Introduction to Analysis of the Infinite Book I},
	  Springer-Verlag, New York, Berlin, Heidelberg, London,
	  Tokyo, 1988. (Orig. 1748)

\vspace{2pt}
\bibitem{LR86}
	 Levesque, Claude; Rhin, Georges.
	\emph{A few classes of periodic continued fractions}. 
		Utilitas Math. 30 (1986), 79--107.

\vspace{2pt}
\bibitem{M01}
	Daniel J. Madden,
	\emph{Constructing Families of Long Continued Fractions},
	Pacific J. Math \textbf{198} (2001), No. 1,
		123--147.
	
\vspace{2pt}
\bibitem{N90}
	Wladyslaw Narkiewicz,
	  \emph{Elementary and Analytic Theory of Algebraic Numbers, 
		(Second Edition)},
	  Springer-Verlag, New York, Berlin, Heidelberg, London,
	  Tokyo, Hong Kong/PWN-Polish Scientific Publishers, Warszawa
1990. 
			(First Edition 1974)


\bibitem{oP13}
   Oskar Perron,
   \emph{Die Lehre von dem Kettenbr\"{u}chen},
   B.G. Teubner, Leipzig-Berlin, 1913.

\vspace{2pt}

\bibitem{S61}
 	Schinzel, A. 
	\emph{On some problems of the arithmetical theory of continued fractions}. 
	Acta Arith. \textbf{6} 1960/1961 393--413.

\vspace{2pt}

\bibitem{S62}
	Schinzel, A. 
	\emph{On some problems of the arithmetical theory of continued fractions. II}. 
	Acta Arith. \textbf{7} 1961/1962 287--298. 

\bibitem{VDP94}
	van der Poorten, A. J. 
	\emph{Explicit formulas for units in certain quadratic number fields.}
	 Algorithmic number theory (Ithaca, NY, 1994), 194--208, 
	Lecture Notes in Comput. Sci., 877, Springer, Berlin, 1994.

\vspace{2pt}



\bibitem{VDPW99}
	van der Poorten, A. J.; Williams, H. C. 
	\emph{On certain continued fraction expansions of fixed period length.}
	 Acta Arith. \textbf{89} (1999), no. 1, 23--35.

\end{thebibliography}
\end{document}